\documentclass{maxart}
\usepackage{amssymb,amsmath,mathrsfs}

\theoremstyle{plain}
\newtheorem{theorem}[equation]{Theorem}
\newtheorem{lemma}[equation]{Lemma}
\newtheorem{corollary}[equation]{Corollary}
\newtheorem{proposition}[equation]{Proposition}

\theoremstyle{definition}
\newtheorem{definition}[equation]{Definition}
\newtheorem{example}[equation]{Example}
\newtheorem{remark}[equation]{Remark}

\newtheorem*{ack}{Acknowledgements}

\theoremstyle{remark}
\newtheorem*{case1}{\indent Case 1}
\newtheorem*{case2}{\indent Case 2}
\newtheorem*{case3}{\indent Case 3}
\newtheorem*{case4}{\indent Case 4}

\newenvironment{proof2}%
	{\begin{proof}}{\renewcommand{\qed}{}\end{proof}}
\newenvironment{enum}%
	{\begin{enumerate}%
	\setlength{\itemsep}{0pt}}%
	{\end{enumerate}}
\newenvironment{pict}[2]%
	{\setlength{\unitlength}{1mm}
	\begin{center}
	\begin{picture}(#1,#2)
	\scriptsize
}%
	{\end{picture}
	\end{center}

	\noindent}

\numberwithin{equation}{section}
\newcommand{\zindex}[3]{\put(#1,#2){\makebox(0,0){${#3}$}}}
\newcommand{\Q}{{\mathbb Q}} 
\newcommand{\R}{{\mathbb R}} 
\newcommand{\Z}{{\mathbb Z}} 
\newcommand{\Gs}{{\mathscr G}}
\newcommand{\abs}[1]{\left\lvert#1\right\rvert}
\DeclareMathOperator{\aut}{Aut}

\renewcommand{\leq}{\leqslant}

\begin{document}

\renewcommand{\baselinestretch}{1.20}
\small

\title[Generalized Baumslag--Solitar groups]{Splittings of
generalized Baumslag--Solitar groups} 
\author{Max Forester}
\address{Mathematics Department, University of Oklahoma, Norman OK 73019,
USA} 
\email{forester@math.ou.edu}

\begin{abstract}
We study the structure of generalized Baumslag--Solitar groups from the
point of view of their (usually non-unique) splittings as fundamental
groups of graphs of infinite cyclic groups. We find and characterize
certain decompositions of smallest complexity (\emph{fully reduced}
decompositions) and give a simplified proof of the existence of
deformations. We also prove a finiteness theorem and solve the
isomorphism problem for generalized Baumslag--Solitar groups with no
non-trivial integral moduli. 
\end{abstract}

\maketitle

\section*{Introduction}

This paper explores the structure of generalized Baumslag--Solitar groups
from the point of view of their (usually non-unique) splittings as
fundamental groups of graphs of groups. By definition, a generalized
Baumslag--Solitar group is the fundamental group of a graph of infinite
cyclic groups. Equivalently, it is a group that acts on a simplicial tree
with infinite cyclic vertex and edge stabilizers. We call such tree
actions generalized Baumslag--Solitar trees. These groups have arisen in
the study of splittings of groups, both in the work of Kropholler
\cite{krop:torus,krop:gbsgroups} and as useful examples of JSJ
decompositions \cite{forester:jsj}. They were classified up to
quasi-isometry in \cite{whyte:gbsgroups,farbmosher:bs1}, but their
group-theoretic classification is still unknown. 

Our approach to understanding generalized Baumslag--Solitar groups is to
study the space of all generalized Baumslag--Solitar trees for a given
group. In most cases this is a \emph{deformation space}, consisting of
$G$-trees related to a given one by a deformation (a sequence of
collapse and expansion moves \cite{herrlich,bass:remarks}). Equivalently
this is the set of $G$-trees having the same elliptic subgroups
(subgroups fixing a vertex) as the given one \cite{forester:trees}. It is
important to note that $G$-trees having the same elliptic subgroups need
not have the same vertex stabilizers. This is one of the main issues
arising in this paper. 

Two notions of complexity for $G$-trees are the number of edge orbits and
the number of vertex orbits. Within a deformation space, the local minima
for both notions occur at the reduced trees: those for which no collapse
moves are possible. In the first part of this paper we study \emph{fully
reduced} $G$-trees. These are $G$-trees in which no vertex stabilizer
contains the stabilizer of a vertex from a different orbit. Fully reduced
trees, when they exist, globally minimize complexity in a deformation
space. They are somewhat canonical (cf. Proposition \ref{VAbijection})
but are not always unique. The slide-inequivalent trees given in
\cite{forester:jsj} are both fully reduced, for example.

One of our main results is Theorem \ref{fullreducibility} which states
that every generalized Baumslag--Solitar tree can be made fully reduced
by a deformation.  After developing properties of fully reduced trees in
Section \ref{sec4} we use these results to give a simplified proof of the
fact (originally proved in \cite{forester:trees}) that all non-elementary
generalized Baumslag--Solitar trees with the same group lie in a single
deformation space.  These results also lay the groundwork for further
study on the classification of generalized Baumslag--Solitar groups. 

In the second part of the paper we focus on generalized Baumslag--Solitar
groups having no non-trivial integral moduli. It turns out that this class
of groups can be understood reasonably well. One key property is given in
Theorem \ref{slidethm}: for such groups, deformations between reduced
trees can be converted into sequences of slide moves, which do not change 
complexity. We then prove a finiteness theorem for such trees (Theorem
\ref{finiteness}), and these two results together yield a solution to the
isomorphism problem (Corollary \ref{isomproblem}). This is our second
main result. 

\begin{ack}
This paper is based on my PhD dissertation, prepared under the direction
of Peter Scott. I wish to express my gratitude to 
Peter Scott for his support and for many
valuable discussions and suggestions. 
I also thank Gilbert Levitt and Noel Brady for helpful discussions and
comments. 
\end{ack}

\section{Preliminaries}\label{sec1}

We will use Serre's notation for graphs and trees
\cite{serre:trees}. Thus a graph $A$ is a pair of sets $(V(A)$, $E(A))$
with maps $\partial_0, \partial_1 \co E(A) \to V(A)$ and an involution $e
\mapsto \overline{e}$ (for $e \in E(A)$), such that
$\partial_i\overline{e} = \partial_{1-i} e$ and $e \not= \overline{e}$
for all $e$. An element $e \in E(A)$ is to be thought of as an oriented
edge with initial vertex $\partial_0 e$ and terminal vertex $\partial_1
e$. We denote by $E_0(v)$ the set of all edges having initial vertex
$v$. An edge $e$ is a \emph{loop} if $\partial_0 e = \partial_1 e$. 

\begin{definition}
Let $G$ be a group. A \emph{$G$-tree} is a tree with a $G$-action by
automorphisms, without inversions.  A $G$-tree is \emph{proper} if every
edge stabilizer is strictly smaller that its neighboring vertex
stabilizers. It is \emph{minimal} if there is no proper $G$-invariant
subtree, and it is \emph{cocompact} if the quotient graph is finite. 
\end{definition}

Given a $G$-tree $X$, an element $g\in G$ is \emph{elliptic} if it fixes
a vertex of $X$ and is \emph{hyperbolic} otherwise. If $g$ is hyperbolic
then there is a unique $G$-invariant line in $X$, called the \emph{axis} of
$g$, on which $g$ acts as a translation. A subgroup $H$ of $G$ is
\emph{elliptic} if it fixes a vertex. 

Suppose a graph of groups has an edge $e$ which is a loop. Let $A$ be the
vertex group and $C$ the edge group, with inclusion maps $i_0, i_1 \co C
\hookrightarrow A$. If one of these maps, say $i_0$, is an isomorphism,
then $e$ is an \emph{ascending loop}. The \emph{monodromy} is the
composition $i_1 \circ i_0^{-1}  \co A \hookrightarrow A$. 

\begin{definition}
In a \emph{collapse move}, an edge in a graph of
groups carrying an amalgamation of the form $A \ast_C C$ is collapsed to
a vertex with group $A$. Every inclusion map having target group $C$ is
reinterpreted as a map into $A$, via the injective map of vertex groups 
$C \hookrightarrow A$. 

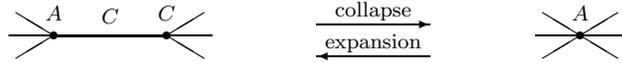
\begin{pict}{90}{10}
\thicklines
\put(10,5){\circle*{1}}
\put(25,5){\circle*{1}}

\put(10,5){\line(1,0){15}}

\thinlines
\put(45,6.5){\vector(1,0){15}}
\put(60,2.2){\vector(-1,0){15}}
\zindex{52.5}{8.2}{\mbox{collapse}}
\zindex{52.5}{4}{\mbox{expansion}}

% \put(25,5){\line(3,5){3}}
\put(25,5){\line(5,-3){5}}
\put(25,5){\line(5,3){5}}

\put(10,5){\line(-5,3){5}}
\put(10,5){\line(-5,-3){5}}

\put(80,5){\circle*{1}}
\put(80,5){\line(-5,3){5}}
\put(80,5){\line(-5,-3){5}}

% \put(80,5){\line(3,5){3}}
\put(80,5){\line(5,-3){5}}
\put(80,5){\line(5,3){5}}

\put(4,5){\line(1,0){27}}
\put(74,5){\line(1,0){12}}

\zindex{10}{8}{A}
\zindex{17.5}{7.5}{C}
\zindex{25}{8}{C}

\zindex{80}{8}{A}
\end{pict}% 
An \emph{expansion move} is the reverse of a collapse move. Both of these
moves are called \emph{elementary moves}. A \emph{deformation} (also
called an \emph{elementary deformation} in
\cite{forester:trees,forester:jsj}) is a finite sequence of such moves. 

A graph of groups is 
\emph{reduced} if it admits no collapse moves. This means that if an
inclusion map from an edge group to a vertex group is an isomorphism,
then the edge is a loop. Correspondingly, a $G$-tree is reduced if,
whenever an edge stabilizer is equal to the stabilizer of one of its
endpoints, both endpoints are in the same orbit. Note that reduced
$G$-trees are minimal. 
\end{definition} 

\begin{definition}
The deformation shown below (cf. \cite{herrlich}) is called a \emph{slide
move}. In order to perform the move it is required that $D \subseteq C$
(regarded as subgroups of $A$). 

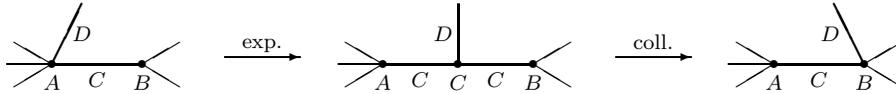
\begin{pict}{120}{12}
\thicklines
\put(102,4){\circle*{1}}
\put(114,4){\circle*{1}}
\put(102,4){\line(1,0){12}}
\put(114,4){\line(-1,2){4}}

\zindex{102}{1.5}{A}
\zindex{114}{1.5}{B}
\zindex{108}{2}{C}
\zindex{109.5}{8}{D}

\thinlines
\put(114,4){\line(5,3){5}}
\put(114,4){\line(5,-3){5}}
\put(102,4){\line(-5,3){5}}
\put(96,4){\line(1,0){6}}
\put(102,4){\line(-5,-3){5}}

%%%%%%%%

\thicklines
\put(6,4){\circle*{1}}
\put(18,4){\circle*{1}}
\put(6,4){\line(1,0){12}}
\put(6,4){\line(1,2){4}}

\zindex{6}{1.5}{A}
\zindex{18}{1.5}{B}
\zindex{12}{2}{C}
\zindex{10}{8}{D}

\thinlines
\put(18,4){\line(5,3){5}}
\put(18,4){\line(5,-3){5}}
\put(6,4){\line(-5,3){5}}
\put(0,4){\line(1,0){6}}
\put(6,4){\line(-5,-3){5}}

%%%%%%%%

\thicklines
\put(50,4){\circle*{1}}
\put(60,4){\circle*{1}}
\put(70,4){\circle*{1}}
\put(50,4){\line(1,0){20}}
\put(60,4){\line(0,1){8}}

\zindex{50}{1.5}{A}
\zindex{55}{2}{C}
\zindex{60}{1.5}{C}
\zindex{65}{2}{C}
\zindex{70}{1.5}{B}
\zindex{58}{8}{D}

\thinlines
\put(70,4){\line(5,3){5}}
\put(70,4){\line(5,-3){5}}
\put(50,4){\line(-5,3){5}}
\put(44,4){\line(1,0){6}}
\put(50,4){\line(-5,-3){5}}

%%%%%%%%

\zindex{34}{6.5}{\mbox{exp.}}
\zindex{86}{7}{\mbox{coll.}}
\put(29,5){\vector(1,0){10}}
\put(81,5){\vector(1,0){10}}

\end{pict}% 
It is permitted for the edge carrying $C$ to be a
loop; in this case the only change to the graph of groups is in the
inclusion map $D \hookrightarrow A$. See Proposition \ref{gbsmoves} for
an example. 
\end{definition}

\begin{definition}
An \emph{induction move} is an expansion and collapse along an ascending
loop. In the diagram below the ascending loop has vertex group $A$ and
monodromy $\phi \co A \to A$, and $B$ is a subgroup such that $\phi(A)
\subseteq B \subseteq A$. The map $\iota\co B \to A$ is inclusion. The
lower edge is expanded and the upper edge is collapsed, resulting in an
ascending loop with monodromy the induced map $\phi\vert_{B}\co B \to B$. 

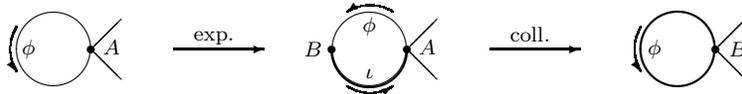
\begin{pict}{104}{12}
\thicklines
\put(92,6){\circle{10}}
\put(97,6){\circle*{1}}

\put(51,6){\oval(10,10)[b]}

\thinlines
\put(9,6){\circle{10}}
\put(14,6){\circle*{1}}

\put(51,6){\circle{10}}
\put(46,6){\circle*{1}}
\put(56,6){\circle*{1}}

\put(14,6){\line(1,1){4}}
\put(14,6){\line(1,-1){4}}

\put(56,6){\line(1,1){4}}
\put(56,6){\line(1,-1){4}}

\put(97,6){\line(1,1){4}}
\put(97,6){\line(1,-1){4}}

\scriptsize
\zindex{5.7}{6}{\phi}
\qbezier(4,3)(2.3,6)(4,9)
\put(4,3){\vector(2,-3){0}}
\zindex{16.7}{6.2}{A}

\put(25,6){\vector(1,0){12}}
\zindex{31}{7.5}{\mbox{exp. }}

\qbezier(48,11)(51,12.7)(54,11)
\put(48,11){\vector(-3,-2){0}}
\zindex{51}{9}{\phi}
\qbezier(48,1)(51,-0.7)(54,1)
\put(54,1){\vector(3,2){0}}
\zindex{51}{2.7}{\iota}
\zindex{43.5}{6}{B}
\zindex{58.7}{6.2}{A}

\put(67,6){\vector(1,0){12}}
\zindex{73}{8}{\mbox{coll. }}

\zindex{88.9}{6}{\phi}
\qbezier(87,3)(85.3,6)(87,9)
\put(87,3){\vector(2,-3){0}}
\zindex{100}{6}{B}

\end{pict}% 
The reverse of this move is also considered an induction move. Notice
that the vertex group changes, in contrast with slide moves. 
\end{definition}

\begin{definition}
A \emph{fold} is most easily described directly in terms of
$G$-trees. The graph of groups description includes many different cases
which are explained in \cite{bestvina:accessibility}. To perform a fold 
in a $G$-tree one chooses edges $e$ and $f$ with $\partial_0 e =
\partial_0 f$, and identifies $e$ and $f$ to a single edge. One also
identifies $g e$ with $g f$ for every $g \in G$, so that
the resulting quotient graph has a $G$-action. It is not difficult to
show that the new graph is a tree. 
\end{definition}

The following basic result is proved in \cite[Proposition
3.16]{forester:trees}. 

\begin{proposition}\label{foldfac}
Suppose a fold between $G$-trees preserves hyperbolicity of
elements of $G$. Then the fold can be represented by a
deformation. \endproof 
\end{proposition}

\section{Generalized Baumslag--Solitar groups}\label{sec2}

\begin{definition}
A \emph{generalized Baumslag--Solitar tree} is a $G$-tree whose vertex
and edge stabilizers are all infinite cyclic. The groups $G$ that arise
are called \emph{generalized Baumslag--Solitar groups}. Basic examples
include Baumslag--Solitar groups \cite{baumslagsolitar}, torus knot and
link groups, and finite index subgroups of these groups. 

The quotient graphs of groups have all vertex and edge groups isomorphic
to $\Z$, and the inclusion maps are multiplication by various non-zero
integers. Thus any example is specified by a graph $A$ and a function
$i\co E(A) \to (\Z - \{0\})$. The corresponding graph of groups will be 
denoted by $(A,i)_{\Z}$. If $X$ is the $G$-tree above $(A,i)_{\Z}$ then
the induced function $i \co E(X) \to (\Z - \{0\})$ satisfies 
\begin{equation}\label{absindex}
\abs{i(e)} \ = \ [G_{\partial_0 e} : G_e] 
\end{equation}
for all $e\in E(X)$. 
\end{definition}

\begin{remark}\label{signs}
There is generally some choice involved in writing down a quotient graph
of groups of a $G$-tree. This issue is explored fully in
\cite[Section 4]{bass:covering}. Without changing the $G$-tree it describes, a
graph of groups may be modified by twisting an inclusion map by an inner
automorphism of the target vertex group. Any two quotient graphs of
groups of a $G$-tree are related by modifications of this type. 

In the case of generalized Baumslag--Solitar trees there are no such
inner automorphisms and the quotient graph of groups is uniquely
determined by the $G$-tree. The associated edge-indexed graph is then very
nearly uniquely determined; the only ambiguity arises from the choice of
generators of edge and vertex groups. One may simultaneously change
the signs of all indices at a vertex, or change the signs of $i(e)$ and
$i(\overline{e})$ together for any $e$, with no change in the graph of
groups or the $G$-tree it encodes. 
\end{remark}

Elementary moves and deformations can be described directly in terms of
edge-indexed graphs, as follows. The verifications are left to the
reader. In the diagrams below, each index $i(e)$ is shown next to the
endpoint $\partial_0 e$. Note in particular that any deformation
performed on a generalized Baumslag--Solitar tree results again in a
generalized Baumslag--Solitar tree. 

\begin{proposition} \label{gbsmoves} 
If an elementary move is performed on a generalized
Baumslag--Solitar tree, then the quotient graph of groups changes locally
as follows: 

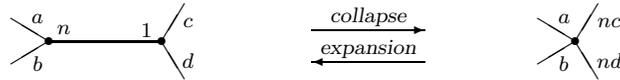
\begin{pict}{90}{10}
\thicklines
\put(10,5){\circle*{1}}
\put(25,5){\circle*{1}}

\put(10,5){\line(1,0){15}}

\thinlines
\put(45,6.5){\vector(1,0){15}}
\put(60,2.2){\vector(-1,0){15}}
\zindex{52.5}{8.2}{\mbox{collapse}}
\zindex{52.5}{4}{\mbox{expansion}}

\put(25,5){\line(3,5){3}}
\put(25,5){\line(3,-5){3}}

\put(10,5){\line(-5,3){5}}
\put(10,5){\line(-5,-3){5}}

\put(80,5){\circle*{1}}
\put(80,5){\line(-5,3){5}}
\put(80,5){\line(-5,-3){5}}

\put(80,5){\line(3,5){3}}
\put(80,5){\line(3,-5){3}}

\scriptsize
\zindex{8.5}{8}{a}
\zindex{8.5}{2}{b}
\zindex{12}{6.5}{n}
\zindex{23}{6.5}{1}
\zindex{28.5}{7.5}{c}
\zindex{28.5}{2.5}{d}

\zindex{78.5}{8}{a}
\zindex{78.5}{2}{b}
\zindex{84.5}{7.5}{nc}
\zindex{84.5}{2.5}{nd}
 
\end{pict} 
A slide move has the following description: 

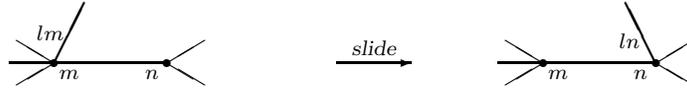
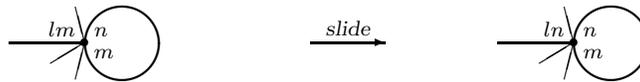
\begin{pict}{100}{11}
\thicklines
\put(75,3){\circle*{1}}
\put(90,3){\circle*{1}}
\put(75,3){\line(1,0){15}}
\put(90,3){\line(-1,2){4}}

\thinlines
\put(47.5,3){\vector(1,0){10}}
\zindex{52.5}{5}{\mbox{slide}}

\put(90,3){\line(5,3){5}}
\put(90,3){\line(5,-3){5}}
\put(75,3){\line(-5,3){5}}
\put(69,3){\line(1,0){6}}
\put(75,3){\line(-5,-3){5}}

\scriptsize
\zindex{77}{1.5}{m}
\zindex{88}{1.5}{n}
\zindex{86.5}{6}{ln}

\thicklines
\put(10,3){\circle*{1}}
\put(25,3){\circle*{1}}
\put(10,3){\line(1,0){15}}
\put(10,3){\line(1,2){4}}

\thinlines
\put(25,3){\line(5,3){5}}
\put(25,3){\line(5,-3){5}}
\put(10,3){\line(-5,3){5}}
\put(4,3){\line(1,0){6}}
\put(10,3){\line(-5,-3){5}}

\scriptsize
\zindex{12}{1.5}{m}
\zindex{23}{1.5}{n}
\zindex{9.5}{7}{lm}
 
\end{pict} 
or

\begin{pict}{100}{10}
\thicklines

\put(82.5,5){\circle*{1}}
\put(87.5,5){\circle{10}}
\put(72.5,5){\line(1,0){10}}

\thinlines
\put(47.5,5){\vector(1,0){10}}
\zindex{52.5}{7}{\mbox{slide}}

\put(82.5,5){\line(-5,-3){4.5}}
\put(82.5,5){\line(-1,-4){1.2}}
\put(82.5,5){\line(-1,4){1.2}}

\scriptsize
\zindex{85.1}{3.5}{m}
\zindex{84.7}{6.5}{n}
\zindex{79.9}{6.8}{ln}

\thicklines

\put(17.5,5){\circle*{1}}
\put(22.5,5){\circle{10}}
\put(7.5,5){\line(1,0){10}}

\thinlines
\put(17.5,5){\line(-5,-3){4.5}}
\put(17.5,5){\line(-1,-4){1.2}}
\put(17.5,5){\line(-1,4){1.2}}

\scriptsize
\zindex{20.1}{3.5}{m}
\zindex{19.7}{6.5}{n}
\zindex{14.5}{6.8}{lm}

\end{pict} 
An induction move is as follows (cf. Lemma \ref{essentialdivides}):

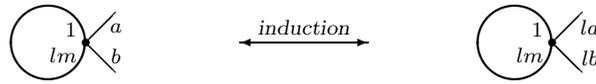
\begin{pict}{80}{10}
\thicklines
\put(6,5){\circle{10}}
\put(11,5){\circle*{1}}

\put(68,5){\circle{10}}
\put(73,5){\circle*{1}}

\thinlines

\put(11,5){\line(1,1){4}}
\put(11,5){\line(1,-1){4}}

\put(73,5){\line(1,1){4}}
\put(73,5){\line(1,-1){4}}

\scriptsize
\zindex{15}{7}{a}
\zindex{15}{3.2}{b}
\zindex{9}{6.7}{1}
\zindex{8}{3.3}{lm}

\put(31.5,5){\vector(1,0){17}}
\put(48.5,5){\vector(-1,0){17}}
\zindex{40}{7}{\mbox{induction}}

\zindex{78}{7}{la}
\zindex{78}{3}{lb}
\zindex{71}{6.7}{1}
\zindex{70}{3.3}{lm}

\end{pict}
\end{proposition}

\begin{definition}
A $G$-tree is \emph{elementary} if there is a $G$-invariant point or
line, and \emph{non-elementary} otherwise. In \cite[Lemma
2.6]{forester:jsj} it is shown that a generalized Baumslag--Solitar 
tree is elementary if and only if the group is isomorphic to $\Z$, $\Z
\times \Z$, or the Klein bottle group. Thus we may speak of generalized
Baumslag--Solitar groups as being elementary or non-elementary. 
\end{definition}

A fundamental property of generalized Baumslag--Solitar groups is that
the elliptic subgroups are canonical, except in the elementary
case. Recall that two subgroups $H, K$ of $G$ are \emph{commensurable} if
$H \cap K$ has finite index in both $H$ and $K$. The following lemma is
proved in \cite[Corollary 6.10]{forester:trees} and \cite[Lemma
2.5]{forester:jsj}. 

\begin{lemma}\label{elliptic}
Let $X$ be a non-elementary generalized Baumslag--Solitar tree with group
$G$. A nontrivial subgroup $H \subseteq G$ is elliptic if and only if it
is infinite cyclic and is commensurable with all of its
conjugates. \endproof 
\end{lemma}

The property of $H$ given in the lemma is rather special. Kropholler
showed in \cite{krop:gbsgroups} that among finitely generated groups of
cohomological dimension $2$, the existence of such a subgroup exactly
characterizes the generalized Baumslag--Solitar groups. 

\section{Full reducibility}\label{sec3}

\begin{definition}
A graph of groups is \emph{fully reduced} if no vertex group can be
conjugated into another vertex group. Correspondingly, a $G$-tree is
fully reduced if, whenever one vertex stabilizer contains another vertex
stabilizer, they are conjugate. Notice that a fully reduced graph of
groups is minimal and reduced. Two basic examples of fully reduced
trees are proper trees and trees having a single vertex orbit. 
\end{definition}

We shall see that for generalized Baumslag--Solitar groups, fully reduced
decompositions exist and have underlying graphs of smallest
complexity (Theorems \ref{fullreducibility} and \ref{complexity} below). 

\begin{example}\label{bs530} 
The $G$-tree shown on the left is reduced but not fully reduced. The
valence three vertex group can be conjugated into the other vertex
group (by conjugating around the loop). After performing an induction
move and a collapse one finds that $G$ is the Baumslag--Solitar group
$BS(5,30)$, which was perhaps not obvious initially. 
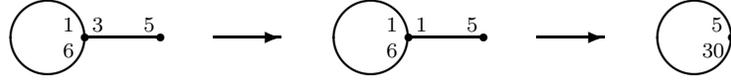
\begin{pict}{100}{10}
\thicklines
\put(7,5){\circle{10}}
\put(12,5){\circle*{1}}
\put(22,5){\circle*{1}}
\put(12,5){\line(1,0){10}}

\put(29,5){\vector(1,0){9}}

\put(50,5){\circle{10}}
\put(55,5){\circle*{1}}
\put(65,5){\circle*{1}}
\put(55,5){\line(1,0){10}}

\put(72,5){\vector(1,0){9}}

\put(93,5){\circle{10}}
\put(98,5){\circle*{1}}

\scriptsize
\zindex{9.8}{6.7}{1}
\zindex{9.8}{3.3}{6}
\zindex{13.7}{6.7}{3}
\zindex{20.5}{6.7}{5}

% \zindex{33.5}{7}{\mbox{ ind. }}

\zindex{52.8}{6.7}{1}
\zindex{52.8}{3.3}{6}
\zindex{56.7}{6.7}{1}
\zindex{63.5}{6.7}{5}

% \zindex{76.5}{7}{\mbox{ coll. }}

\zindex{96}{6.7}{5}
\zindex{95.5}{3.3}{30}

\end{pict}%
The following result generalizes this procedure to arbitrary generalized
Baumslag--Solitar trees. 
\end{example}

\begin{theorem}\label{fullreducibility}
Every cocompact generalized Baumslag--Solitar tree is related by a
deformation to a fully reduced (generalized Baumslag--Solitar) tree. 
\end{theorem}

The proof relies strongly on the fact that stabilizers are infinite
cyclic, and therefore contain a unique subgroup of any given
index. Before proving the theorem we establish some preliminary facts
concerning generalized Baumslag--Solitar trees. Our first objective
(Corollary~\ref{fixpathcor}) is to
characterize the paths that are fixed by vertex stabilizers. 

\begin{lemma}\label{fixpath}
Let $X$ be a generalized Baumslag--Solitar tree with group $G$. Suppose $G_x
\subseteq n G_{x'}$ for vertices $x \not= x'$. Let $(e_1, \ldots, e_k)$ be
the path from $x$ to $x'$, with vertices $x_0 = x$, $x_i = \partial_1
e_i$ for $1 \leq i \leq k$. Define $m_i = i(\overline{e}_i)$, $n_i =
i(e_{i+1})$, and $n_k = n$. Then $i(e_1) = \pm 1$ and
\begin{enumerate}
\item[\textup{(i)}] $G_x = \left(\Pi_{i=1}^{k} \, m_i \, / \, \Pi_{i=1}^{k-1}
\,  n_i\right)   G_{x_k}$  
\item[\textup{(ii)}] $\Pi_{i=1}^k \ n_i$ divides $\Pi_{i=1}^k \ m_i$. 
\end{enumerate}
\end{lemma}

\begin{proof}
The statement $i(e_1)=\pm 1$ is clear because $G_x$ fixes $e_1$. 
The other two statements are proved together by induction on $k$. 

If $k=1$ then (i) says that $G_x = m_1 G_{x'}$, which holds because
$i(e_1) = \pm 1$. Then the assumption $G_x \subseteq n_1 G_{x'}$ implies
that $n_1$ divides $m_1$, because $G_x = G_{e_1} = m_1 G_{x'}$. 

Now let $k > 1$ be arbitrary. Since $G_x$ fixes $e_k$ we have $G_x
\subseteq n_{k-1} G_{x_{k-1}}$, and the induction hypothesis gives that
$\Pi_{i=1}^{k-1} n_i$ divides $\Pi_{i=1}^{k-1} m_i$. We also have 
$G_x = \left(\Pi_{i=1}^{k-1} \, m_i \, / \, \Pi_{i=1}^{k-2}
\,  n_i\right) 
G_{x_{k-1}}$ and $n_{k-1} G_{x_{k-1}} = G_{e_k} = m_k
G_{x_k}$. Therefore \[G_x \ = \ \left(\Pi_{i=1}^{k-1} \, m_i \, / \,
\Pi_{i=1}^{k-1} \,  n_i\right) n_{k-1} \, G_{x_{k-1}} \ = \
\left(\Pi_{i=1}^{k-1} \, m_i \, / \,
\Pi_{i=1}^{k-1} \,  n_i\right) m_k \, G_{x_k},\] proving (i). 

Next, the assumption $G_x \subseteq n_k G_{x_k}$ becomes
$\left(\Pi_{i=1}^{k} \, m_i \, / \, \Pi_{i=1}^{k-1} \,  n_i\right)
G_{x_k} \subseteq n_k G_{x_k}$ by (i), establishing (ii). 
\end{proof}

\begin{remark} 
The lemma is valid in any locally finite $G$-tree,
provided one interprets statements such as $G_x \subseteq n G_{x'}$
correctly. For example this statement would mean that $G_x \subseteq G_{x'}$
and $n$ divides $[G_{x'}:G_x]$. The following corollary, however, is
specific to generalized Baumslag--Solitar trees. 
\end{remark}

\begin{corollary}\label{fixpathcor}
Let $(e_1, \ldots, e_k)$ be a path in $X$ and define $m_i$, $n_i$ as in
the previous lemma. Then ${G}_{\partial_0 e_1}$ fixes the path $(e_1,
\ldots, e_k)$ if and only if $i(e_1) = \pm 1$ and for every $l\leq
(k-1)$ 
\begin{equation}\label{fixpathineq2} 
\Pi_{i=1}^l \ n_i \quad \text{divides} \quad \Pi_{i=1}^l \ m_i. 
\end{equation}
\end{corollary}

\begin{proof}
The forward implication is given by Lemma \ref{fixpath}(ii). The converse
is proved by induction on $k$. Suppose \eqref{fixpathineq2} holds
for each $l$ and that ${G}_{\partial_0 e_1}$ fixes the path $(e_1,
\ldots, e_{k-1})$. Then $G_{\partial_0 e_1}$ is the subgroup of
$G_{\partial_0 e_k}$ of index $\Pi_{i=1}^{k-1} m_i \, / \,
\Pi_{i=1}^{k-2}n_i$  by Lemma~\ref{fixpath}(i). Property
\eqref{fixpathineq2} for $l= k-1$ implies that $n_{k-1}$ divides this
index, and so ${G}_{\partial_0 e_1} \subseteq G_{e_k}$. 
\end{proof}

Next we describe the steps needed to construct the deformation of
Theorem~\ref{fullreducibility}. The kind of example one should have in mind
is one similar to Example \ref{bs530}, but with several loops incident to
the left-hand vertex. 

Throughout the rest of this section $X$ denotes a generalized
Baumslag--Solitar tree with group $G$ and quotient graph of groups
$(A,i)_{\Z}$. 

\begin{definition}
Let $f\in E(A)$ be an edge with $\partial_0 f \not= \partial_1 f$. Let
$\rho = (e_1, \ldots, e_k, f)$ be a path in $A$ such that each $e_i$ is a
loop at $\partial_0 f$. We say that $\rho$ is an \emph{admissible path
for $f$} if, for some lift $\tilde{\rho} = 
(\tilde{e}_1, \ldots, \tilde{e}_k, \tilde{f}) \subset X$,
\[{G}_{\partial_0 \tilde{e}_1} \ \ \subseteq \ \
{G}_{\tilde{f}}.\] 
This condition depends only on the indices along the path $\rho$ by
Corollary~\ref{fixpathcor}, so it is independent of the choice of
$\tilde{\rho}$. When dealing with admissible paths we will use
the notation $m_i = i(\overline{e}_i)$, $n_i = i(e_{i+1})$, and $n_k =
i(f)$; then the path is admissible if and only if $i(e_1) = \pm 1$ and
\eqref{fixpathineq2} holds for each $l$. 

An edge $e\in E(A)$ is \emph{essential} if $i(e) \not= \pm 1$, and
\emph{inessential} otherwise. The \emph{length} of $\rho$ is $k$, and the
\emph{essential length} of $\rho$ is the number of essential edges
occurring in $\rho$. 
\end{definition}

\begin{lemma}\label{inessentialfirst}
If $\rho = (e_1, \ldots, e_k, f)$ is an admissible path then there is a 
permutation $\sigma$ such that the path $\rho_{\sigma} = (e_{\sigma(1)},
\ldots, e_{\sigma(k)}, f)$ is admissible and all of the essential
edges of $\rho_{\sigma}$ occur after the inessential edges.
\end{lemma}

\begin{proof}
We show first that if $e_j$ is inessential for some $j > 1$ then 
\[\rho' = \rho_{((j-1) \ j)} = (e_1, \ldots, e_{j-2}, e_j, e_{j-1}, e_{j+1},
\ldots, e_k, f)\] is an admissible path for $f$. Letting $m'_i$ and
$n'_i$ be the indices along $\rho'$, we have $m'_{j-1} = m_j$, $m'_j =
m_{j-1}$, and $n'_{j-2} = n_{j-1}$, $n'_{j-1} = n_{j-2}$, with all other
indices unchanged. Clearly \eqref{fixpathineq2} still holds for $l \not=
j-2, j-1$. One easily verifies \eqref{fixpathineq2} for these other two
cases as well, using the fact that $n_{j-1} = \pm 1$ (because $e_j$ is
inessential). Hence $\rho'$ is admissible for $f$. Next, by
using transpositions of this type, one can move all of the inessential
edges in $\rho$ to the front of the path. 
\end{proof}

\begin{lemma}\label{inessentialsame}
Let $\rho = (e_1, \ldots, e_k, f)$ be an admissible path such that $e_1,
\ldots, e_j$ are inessential and $e_{j+1}, \ldots, e_k$ are
essential. Then there is a sequence of slide moves, after which the
path $\rho' = (e_1, \ldots, e_1, e_{j+1}, \ldots, e_k, f)$ is
admissible. The inessential part $(e_1, \ldots, e_1)$ of $\rho'$ may have
length greater than $j$, though $\rho$ and $\rho'$ have the same
essential length ($k-j$). 
\end{lemma}

\begin{proof}
First we slide $\overline{e}_1$ over each edge of $(A,i)_{\Z}$ (other
than $e_1$) that appears in $(e_2, \ldots, e_j)$. Since these edges are all
inessential loops, these slides can be performed. The index $i(e_1)$ is
unchanged so occurrences of $e_1$ in the path are still inessential. 

Each slide of $\overline{e}_1$ over $e_i$ multiplies $i(\overline{e}_1)$
by $\pm i(\overline{e}_i)$. The end result is that $m_1$ gets multiplied
by a product $\pm \Pi_{\nu=1}^r n_{i_{\nu}}$. By itself this change does
not violate the conditions \eqref{fixpathineq2}. However if some $e_i$ is
equal to $\overline{e}_1$ (where $i > j$), then $n_{i-1}$ is also
multiplied by $\pm \Pi_{\nu=1}^r n_{i_{\nu}}$, and \eqref{fixpathineq2}
may fail. To remedy this we adjoin several copies of $e_1$ to
the front of the path, one for each occurrence of $\overline{e}_1$ in
$(e_{j+1}, \ldots, e_k)$. Then the products $\Pi_{i=1}^l m_i$ acquire
enough additional factors $\pm \Pi_{\nu=1}^r n_{i_{\nu}}$ to remain
divisible by $\Pi_{i=1}^l n_i$. This extended path is therefore
admissible. As a result of the slide moves, $m_i$ now divides $m_1$
for each $i \leq j$. 

We now replace $e_i$ by $e_1$ for $i \leq j$. The quantities
$\Pi_{i=1}^l m_i$ increase and each $\Pi_{i=1}^l n_i$ remains unchanged
(up to sign), so property \eqref{fixpathineq2} still holds for every~$l$. 
\end{proof}

\begin{lemma}\label{essentialdivides}
Let $(e_1, \ldots, e_1, e_{j+1}, \ldots, e_k, f)$ be an admissible path
such that $e_1$ is inessential and $e_{j+1}, \ldots, e_k$ are
essential. Then there is a sequence of induction moves after which a path
of the same form is admissible, has essential length at most ($k -
j$), and satisfies $i(e_{j+1}) = \pm i(\overline{e}_1)^r$ for some $r$. 
\end{lemma}

\begin{proof}
If $e_{j+1} = \overline{e}_1$ then we can discard it from the path
without affecting admissibility. Thus we can assume that $e_{j+1} \not=
\overline{e}_1$. 
Admissibility implies that $i(e_{j+1})$ divides $i(\overline{e}_1)^j$. 
Let $r$ be minimal so that $i(e_{j+1})$
divides $i(\overline{e}_1)^r$ and let $l$ be any factor of
$i(\overline{e}_1)^r / i(e_{j+1})$ that divides $i(\overline{e}_1)$. We
show how to make $i(e_{j+1})$ become $l \cdot i(e_{j+1})$. By repeating 
this procedure the desired result can be achieved.

Writing $i(\overline{e}_1)$ as $lm$, we perform an induction move along
$e_1$ as follows: 
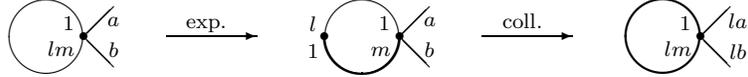
\begin{pict}{100}{10}
\thicklines
\put(48,5){\oval(10,10)[b]}

\put(88,5){\circle{10}}
\put(93,5){\circle*{1}}

\thinlines
\put(6,5){\circle{10}}
\put(11,5){\circle*{1}}

\put(48,5){\circle{10}}
\put(43,5){\circle*{1}}
\put(53,5){\circle*{1}}

\put(11,5){\line(1,1){4}}
\put(11,5){\line(1,-1){4}}

\put(53,5){\line(1,1){4}}
\put(53,5){\line(1,-1){4}}

\put(93,5){\line(1,1){4}}
\put(93,5){\line(1,-1){4}}

\scriptsize
\zindex{15}{7}{a}
\zindex{15}{3.2}{b}
% \zindex{-1}{4.8}{e_1}
\zindex{9}{6.7}{1}
\zindex{8}{3.3}{lm}

\put(22,5){\vector(1,0){12}}
\zindex{28}{6.5}{\mbox{exp. }}% e'_1}

\zindex{57}{7}{a}
\zindex{57}{3.2}{b}
% \zindex{42}{9}{e_1}
% \zindex{42}{1}{e'_1}
\zindex{51}{6.7}{1}
\zindex{50.5}{3}{m}
% \zindex{45.5}{6.7}{l}
% \zindex{45.5}{3.3}{1}
\zindex{41.5}{6.7}{l}
\zindex{41.5}{3}{1}

\put(64,5){\vector(1,0){12}}
\zindex{70}{7}{\mbox{coll. }}% e_1}

\zindex{98}{7}{la}
\zindex{98}{3}{lb}
% \zindex{81}{5}{e'_1}
\zindex{91}{6.7}{1}
\zindex{90}{3.3}{lm}

\end{pict}%
The index of every edge incident to $\partial_0 e_1$ is multiplied by
$l$, except for $i(e_1)$ and $i(\overline{e}_1)$, which remain the
same. As a result the indices $m_i$ and $n_{i-1}$ are multiplied by $l$
whenever $e_i$ is not equal to $e_1$ or $\overline{e}_1$. For every such
$i$ we adjoin a copy of $e_1$ to the front of the path, making it
admissible as in the proof of the preceding lemma. 
\end{proof}

\begin{remark}\label{esslengthone}
The previous argument is still valid when $j=k$. That is, if the path $(e_1,
\ldots, e_1, f)$ is admissible, then there is a sequence of induction
moves after which $i(f) = i(\overline{e}_1)^r$ for some $r$. 
\end{remark}

\begin{proof}[Proof of Theorem~\ref{fullreducibility}.]
We show that if $(A,i)_{\Z}$ is not fully reduced then there is a
deformation to a decomposition having fewer edges. Repeating the
procedure will eventually produce a fully reduced decomposition. 

If $(A,i)_{\Z}$ is not fully reduced then there exist vertices $v$,
$w$ of $X$ such that ${G}_v \subseteq {G}_w$ and $v \not\in G w$. The path
$(\tilde{e}_1, \ldots, \tilde{e}_r)$ from 
$v$ to $w$ contains an edge that does not map to a loop in $A$. Let
$\tilde{f} = \tilde{e}_{k+1}$ be the first such edge. Since ${G}_v$
stabilizes $(\tilde{e}_1, \ldots, \tilde{e}_k, \tilde{f})$, the image 
$\rho = (e_1, \ldots, e_k, f)$ of this path in $A$ is an admissible path
for $f$.  

Next we show how to produce an admissible path for $f$ having essential
length smaller than that of $\rho$, assuming this length is positive. 
Suppose $\rho$ has essential length $s$. Applying
Lemmas~\ref{inessentialfirst},~\ref{inessentialsame},
and~\ref{essentialdivides} 
in succession to the path $\rho$ we can arrange that there is an
admissible path 
$\rho' = (e'_1, \ldots, e'_1, e'_{k'-s+1}, \ldots, e'_{k'}, f)$ 
such that $e'_1$ is inessential, $e'_{k'-s+1}, \ldots, e'_{k'}$ are
essential, and $i(e'_{k'-s+1}) = i(\overline{e}'_1)^r$ for some $r$. 
(The essential length $s$ may have decreased, but then we are done for
the moment.) In applying these lemmas the decomposition $(A,i)_{\Z}$
changes by a deformation to $(A', i')_{\Z}$ where $A'$ has the same
number of edges as $A$. Now we can slide $e'_{k'-s+1}$ over
$\overline{e}'_1$ $r$ times to make $i(e'_{k'-s+1}) = \pm 1$. 

These slide moves affect the indices of the edges $e'_i$ that are equal
to $e'_{k'-s+1}$ or $\overline{e}'_{k'-s+1}$. If $e'_i = e'_{k'-s+1}$
then $n_{i-1}$ is divided by $(m_1)^r = i(\overline{e}'_1)^r$ and this
change does not affect the admissibility of $\rho'$. If $e'_i =
\overline{e}'_{k'-s+1}$ then $m_i$ is divided by $(m_1)^r$. In order to
keep $\rho'$ admissible we adjoin $r$ copies of $e'_1$ to the front of
the path for each such $e'_i$. This
done, we have produced an admissible path for $f$ of smaller essential
length because $e'_{k'-s+1}$ is now inessential. 

By repeating this process we can obtain a decomposition $(A',i')_{\Z}$
related by a deformation to $(A,i)_{\Z}$ (with $\abs{E(A')} =
\abs{E(A)}$), and an 
admissible path for $f$ having essential length zero. Applying
Lemmas~\ref{inessentialsame} and~\ref{essentialdivides} once more, this
path has the form $(e'_1, \ldots, e'_1, f)$ where $i(f) =
i(\overline{e}'_1)^r$. Now we slide $f$ over $\overline{e}'_1$ $r$ times
to make $i(f) = \pm 1$, and collapse $f$. The resulting decomposition 
has fewer edges than $(A,i)_{\Z}$. 
\end{proof}

\section{Vertical subgroups}\label{sec4}

In this section we link the structure of a fully reduced $G$-tree to that
of the group $G$, using the notion of a vertical subgroup. We are
concerned with the difference between elliptic subgroups (which may be
uniquely determined) and vertex stabilizers (which often are not). It
turns out that vertical subgroups are a useful intermediate notion. See in 
particular Example \ref{verticaleg}. 

\begin{definition}
Let $X$ be a $G$-tree. 
A subgroup $H \subseteq G$ is \emph{vertical} if it is
elliptic and every elliptic subgroup containing $H$ is
conjugate to a subgroup of $H$. 
\end{definition}

\begin{lemma}\label{vertical}
If $X$ is fully reduced then an elliptic subgroup is 
vertical if and only if it contains a vertex stabilizer. 
\end{lemma}

\begin{proof}
Suppose $H$ contains $G_v$, and let $H'$ be an elliptic subgroup
containing $H$. Since $H'$ is elliptic, it is contained in $G_w$ for some
$w$, and hence $G_v \subseteq G_w$. Full reducibility implies that $(G_w)^g
= G_v$ for some $g$, and therefore $(H')^g \subseteq H$. 

Conversely suppose $H$ is vertical. Then $H\subseteq G_v$ for some $v$, and
so $(G_v)^g \subseteq H$ for some $g\in G$. Hence $H$ contains
$G_{gv}$. 
\end{proof}

\begin{example}\label{verticaleg}
Let $G$ be the Baumslag--Solitar group $BS(1,6)$ with its standard
decomposition $G = \Z \,\ast_{\, \Z}$ and presentation $\langle x, t
\mid t x t^{-1} = x^6 \rangle$. The vertex stabilizers of the Bass--Serre
tree $X$ are the conjugates of the subgroup $\langle x \rangle$. Among
the subgroups of the form $\langle x^n \rangle$, notice that all are
elliptic, and only those where $n$ is a power of $6$ are vertex
stabilizers. According to Lemma \ref{vertical}, $\langle x^n \rangle$
is vertical if and only if $n$ divides a power of $6$. 

Now consider the automorphism $\phi\co {G} \rightarrow {G}$ defined by
$\phi(x) = x^3$, $\phi(t) = t$ (with inverse $x \mapsto t^{-1} x^2 t$,
$t\mapsto t$). If we twist the action of $G$ on $X$ by $\phi$
then the vertex stabilizers will be the conjugates of $\langle x^2
\rangle$ rather than $\langle x \rangle$. Thus there is no hope of
characterizing vertex stabilizers from the structure of $G$ alone. On
the other hand, the vertical subgroups are uniquely determined (because
the elliptic subgroups are). 

In this particular example, the set of vertical subgroups is the
smallest $\aut(G)$-invariant set of elliptic subgroups containing a
vertex stabilizer. Every vertical subgroup can be realized as a
vertex stabilizer by twisting by an automorphism. 
\end{example}

\begin{definition}\label{eqreln1}
Now we define an equivalence relation on the set of vertical subgroups of
$G$. Set $H \sim K$ if $H$ is conjugate to a subgroup of $K$. 
This relation is symmetric: suppose $H \sim K$, so that $H^g \subseteq K$
for some $g\in G$. Then $H \subseteq K^{g^{-1}}$ and so $K^{g^{-1}}$ is
conjugate to a subgroup of $H$, as $H$ is vertical. Therefore $K \sim
H$. Reflexivity and transitivity are clear.  
\end{definition}

\begin{proposition}\label{VAbijection}
If $X$ is fully reduced then the vertex orbits correspond bijectively
with the $\sim$-equivalence classes of vertical subgroups of $G$. The
bijection is induced by the natural map $v \mapsto G_v$. 
\end{proposition}

\begin{proof}
The induced map is well defined since $G_v \sim (G_v)^g = G_{gv}$ for any
$g$. For injectivity, suppose that $G_v \sim G_w$. Then $G_{gv} \subseteq
G_w$ for some $g\in G$. Full reducibility implies that $gv$ and $w$ are
in the same orbit, hence $v$ and $w$ are as well. 

For surjectivity, suppose $H$ is vertical. It contains a stabilizer
$G_v$ by Lemma \ref{vertical}, and $G_v \subseteq H$ implies $G_v \sim
H$. 
\end{proof}

The following application of Proposition \ref{VAbijection} explains the
choice of the term \emph{fully reduced}. 

\begin{theorem}\label{complexity}
A non-elementary cocompact generalized Baumslag--Solitar tree is fully
reduced if and only if it has the smallest number of edge orbits among
all generalized Baumslag--Solitar trees having the same group. 
\end{theorem}

\begin{proof}
The proof of Theorem \ref{fullreducibility} shows that any generalized
Baumslag--Solitar tree with the smallest number of edge orbits is fully
reduced. For the converse we show that no tree with more edge orbits can
also be fully reduced. 

Suppose the given tree is fully reduced. Let $N\subseteq G$ be the normal
closure of the set of elliptic elements. This subgroup is uniquely
determined since the tree is non-elementary. Note that $G/N$ is the
fundamental group of the quotient graph. Hence the homotopy type of this
graph is uniquely determined. Proposition \ref{VAbijection} implies that
the number of vertices is also uniquely determined, and so the number of
edges is as well. Thus any two fully reduced trees have the same
number of edge orbits. 
\end{proof}

\section{Existence of deformations} \label{sec5}

We now know that generalized Baumslag--Solitar trees can be made fully
reduced (Theorem \ref{fullreducibility}) and that for such trees, the
structure of the tree is partially encoded in the set of elliptic
subgroups (Proposition \ref{VAbijection}). Using these facts we may now
give a quick proof of the existence of deformations between generalized
Baumslag--Solitar trees.  This result is a special case of 
Theorem 1.1 of \cite{forester:trees}. 

\begin{theorem}\label{defthm}
Let $X$ and $Y$ be non-elementary cocompact generalized Baumslag--Solitar
trees with isomorphic groups. Then $X$ and $Y$ are related by a
deformation. 
\end{theorem}

\begin{definition} 
A map between trees is a \emph{morphism} if it sends vertices to vertices
and edges to edges (and respects the maps $\partial_0$, $\partial_1$, $e 
\mapsto \overline{e}$). Geometrically it is a simplicial map which does
not send any edge into a vertex. 
\end{definition}

The following result is taken from \cite[Section 2]{bestvina:accessibility}. 

\begin{proposition}[Bestvina--Feighn]\label{bfprop} 
Let $G$ be a finitely generated group and suppose that $\phi\co X
\rightarrow Y$ is an equivariant morphism of $G$-trees. Assume
further that $X$ is cocompact, $Y$ is minimal, and the edge stabilizers
of $Y$ are finitely generated. Then $\phi$ is a finite composition of 
folds. 
\endproof 
\end{proposition}

\begin{proof}[Proof of Theorem \ref{defthm}.] 
Let $G$ be the common group acting on $X$ and $Y$. Note that both trees
define the same elliptic and vertical subgroups. By Theorem
\ref{fullreducibility} we can assume that both trees are fully reduced (and
minimal). Applying Propositions \ref{bfprop} and \ref{foldfac}, it now
suffices to construct a morphism from $X$ to $Y$. In fact we shall
construct such a map from $X'$ to $Y$, where $X'$ is obtained from $X$ by
subdivision (a special case of a deformation). 

Let $x_1, \ldots,  x_n \in V(X)$ be representatives of the vertex orbits
of $X$. Then there are vertices $y_1, \ldots, y_n \in
V(Y)$ such that $G_{x_i} \subseteq G_{y_i}$, since each $G_{x_i}$ is
elliptic. We define a map $\phi\co V(X) \to V(Y)$ by setting $\phi(x_i)
= y_i$ and extending equivariantly. We then extend $\phi$ to a
topological map $X \to Y$ by sending an edge $e$ to the unique reduced
path in $Y$ from $\phi(\partial_0 e)$ to $\phi(\partial_1
e)$. Subdividing where necessary, we obtain an equivariant simplicial map 
$\phi' \co X' \to Y$. 

Now we verify that $\phi'$ is a morphism. It suffices to check that
$\phi(x) \not= \phi(x')$ whenever $x$ and $x'$ are vertices of $X$ that
bound an edge. There are two cases. If $x' = gx$ for some $g\in G$ then
$g$ is hyperbolic, since it has translation length one in $X$, and
equivariance implies that $\phi(x) \not= \phi(x')$. Otherwise, if $x$ and
$x'$ are in different orbits, then $G_x \not\sim G_{x'}$ by Proposition
\ref{VAbijection}. Here we are using the fact that $X$ is fully reduced.
Equivariance yields $G_{x} \subseteq G_{\phi(x)}$ and $G_{x'} \subseteq
G_{\phi(x')}$, and since these are all vertical subgroups (by Lemma
\ref{vertical}) we now have $G_{x} \sim G_{\phi(x)}$ and $G_{x'} \sim
G_{\phi(x')}$. Hence $G_{\phi(x)} \not\sim G_{\phi(x')}$, and in
particular $\phi(x) \not= \phi(x')$. 
\end{proof}

\section{The modular homomorphism}\label{sec6}

Let ${\Q}^{\times}_{>0}$ denote the positive rationals considered as a
group under multiplication. The following notion was first defined by
Bass and Kulkarni \cite{bass:treelat}. 

\begin{definition}
The \emph{modular homomorphism} $q\co G \rightarrow {\Q}^{\times}_{>0}$
of a locally finite $G$-tree is given by 
\[q(g) \ = \ [V:V\cap V^{g}] \ / \ [V^{g}: V \cap
V^{g}] \] 
where $V$ is any subgroup of $G$ commensurable with a vertex
stabilizer. In this definition we are using the fact that in locally
finite $G$-trees, vertex stabilizers are commensurable with all of
their conjugates. One can easily check that $q$ is independent of the
choice of $V$. 
\end{definition}

In the case of generalized Baumslag--Solitar trees the modular
homomorphism may be defined directly in terms of the graph of groups
$(A,i)_{\Z}$, as in \cite{bass:treelat}. First note that $q$ factors
through $H_1(A)$ because it is trivial on elliptic subgroups and
$\Q^{\times}_{>0}$ is abelian. Writing $q$ as a composition $G \to
H_1(A) \to \Q^{\times}_{>0}$, the latter map is then given by 
\begin{equation}\label{modhom} 
(e_1, \ldots, e_k) \ \mapsto \ \Pi_{j=1}^k \, \abs{i(e_j) / 
i(\overline{e}_j)} .
\end{equation}
To verify \eqref{modhom} note that given $g\in G$, the corresponding
$1$-cycle in $H_1(A)$ is obtained by projecting any (oriented) segment of
the form $[v,gv]$ to $A$. One then uses $V = G_v$ to evaluate $q(g)$,
by applying \eqref{absindex} to the edges of $[v,gv]$. 

The next definition is not actually needed in this paper. We mention it
for completeness, with the expectation that it will be useful in future
work. 

\begin{definition}\label{signedmodhom}
The \emph{signed modular homomorphism} $\hat{q} \co G \to \Q^{\times}$ 
of a generalized Baum\-slag--Solitar tree with quotient graph of groups
$(A,i)_{\Z}$ is defined via the map $H_1(A) \to \Q^{\times}$ given by 
\begin{equation}\label{smodhom} 
(e_1, \ldots, e_k) \ \mapsto \ \Pi_{j=1}^k \, i(e_j) /
i(\overline{e}_j) .
\end{equation}
One should verify that this is well defined in light of Remark
\ref{signs}. Clearly, changing the signs of $i(e)$ and $i(\overline{e})$
together, for any $e$, has no effect. Similarly, since $(e_1, \ldots,
e_k)$ is a cycle, changing all signs at a vertex will introduce an even
number of sign changes in \eqref{smodhom}. 

There is also an \emph{orientation homomorphism} $G \to \{\pm 1\}$ defined
by $g \mapsto \hat{q}(g)/q(g)$. 
\end{definition}

\begin{remark}\label{modhominvariant}
The modular homomorphisms are invariant under deformations. For the
unsigned case, note that that during an elementary move there is a vertex
stabilizer that remains unchanged. Taking $V$ to be this stabilizer, one
obtains invariance of $q$. Alternatively, one may verify directly that
the homomorphisms defined by \eqref{modhom} and \eqref{smodhom} are
invariant, using Proposition \ref{gbsmoves}.  
\end{remark}

\section{Deformations and slide moves} \label{sec7}

In this section we show how to rearrange elementary moves between
generalized Baumslag--Solitar trees.  Our goal is to replace deformations
by sequences of slide moves, which are considerably easier to work
with. It should be noted that in general, reduced generalized
Baumslag--Solitar trees with the same group $G$ need not be related by
slide moves; see \cite{forester:jsj}. Nevertheless this does occur in a
special case, given in Theorem \ref{slidethm} below. 

\begin{definition}
Suppose $(A,i)_{\Z}$ has a loop $e$ with $(i(e), i(\overline{e})) =
(m,n)$. If $m$ divides $n$ then $e$ is a \emph{virtually ascending loop}. It
is \emph{strict} if $n \not= \pm m$. Similarly, a \emph{strict ascending
loop} is one with indices of the form $(\pm 1, n)$, $n \not= \pm 1$. 
\end{definition}

The next two propositions are valid for sequences of moves between
generalized Baumslag--Solitar trees. 

\begin{proposition}\label{expslide}
Suppose an expansion is followed by a slide move. Either
\begin{enum}
\item[\textup{(i)}] the moves remove a strict virtually ascending loop and
create a strict ascending loop, or 
\item[\textup{(ii)}] the moves may be replaced by a (possibly empty)
sequence of slides, followed by an expansion. 
\end{enum}
\end{proposition}

\begin{proof2}
Suppose the expansion creates $e$ and the second move slides
$e_0$ over $e_1$ (from $\partial_0 e_1$ to $\partial_1 e_1$). If $e$ is
not $e_i$ or $\overline{e}_i$ ($i=0,1$) then the moves may be performed
in reverse order as they do not interfere with each other. 

If $e = e_1$ or $\overline{e}_1$ then there is no need to perform the
slide at all.  When performing an expansion at a vertex, the incident
edges are partitioned into two sets, which are then separated by a new
edge. Sliding $e_0$ over the newly created edge is equivalent to
including $e_0$ in the other side of the partition before expanding. 

If $e=e_0$ or $\overline{e}_0$ then there are several cases to
consider. Orient $e$ so that $i(e)=1$, and $e_1$ so that $e$ slides over
$e_1$ from $\partial_0(e_1)$ to $\partial_1(e_1)$. The cases depend on
which of the vertices $\partial_0 e$, $\partial_1 e$, $\partial_0 e_1$,
$\partial_1 e_1$ coincide (after the expansion and before the slide). 

\begin{case1} After expanding $e$, $e_1$ has endpoints $\partial_0(e)$ and
$\partial_1(e)$. If $\partial_1(e) = \partial_0(e_1)$ and $\partial_0(e)
= \partial_1(e_1)$ then $i(e)$ is
still $1$ after the slide and $e$ has become an ascending loop. In
addition, since the slide takes place we must have $i(e_1) \mid
i(\overline{e})$. Writing $i(e_1) = k$ and $i(\overline{e}) = kl$, we
must have had $kl \mid i(\overline{e}_1)$ and $i(e_1) = k$ before the
expansion, so $e_1$ was a virtually ascending loop. Writing $i(\overline{e}_1)
= klm$ (before the expansion) the modulus of the loop $e_1$ is $lm$. If
$lm \not= \pm 1$ then alternative (i) holds. If $lm = \pm 1$ then after
the expansion, $i(\overline{e}_1) = \pm 1$. The expansion and slide may
then be replaced by a single expansion. 

Otherwise
$\partial_1(e) = \partial_1(e_1)$ and $\partial_0(e) =
\partial_0(e_1)$. The two moves have the form:
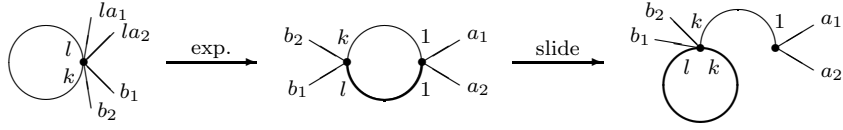
\begin{pict}{111}{16}
\scriptsize
\thicklines
\put(10,8){\circle*{1}}
\put(45,8){\circle*{1}}
\put(55,8){\circle*{1}}
\put(50,8){\oval(10,10)[b]}

\put(92,5){\circle{10}}
\put(92,10){\circle*{1}}
\put(102,10){\circle*{1}}

\thinlines
\put(5,8){\circle{10}}
\put(50,8){\circle{10}}
\put(97,10){\oval(10,10)[t]}

\put(21,8){\vector(1,0){12}}
\put(67,8){\vector(1,0){12}}
\zindex{27}{9.5}{\mbox{exp.}}
\zindex{73}{10}{\mbox{slide}}
\put(10,8){\line(1,6){1}}
\put(10,8){\line(1,1){4}}
\put(10,8){\line(1,-1){4}}
\put(10,8){\line(1,-6){1}}

\put(45,8){\line(-5,3){5}}
\put(45,8){\line(-5,-3){5}}
\put(55,8){\line(5,3){5}}
\put(55,8){\line(5,-3){5}}

\put(92,10){\line(-1,1){4}}
\put(92,10){\line(-6,1){6}}
\put(102,10){\line(5,3){5}}
\put(102,10){\line(5,-3){5}}

\zindex{44.5}{11.5}{k}
\zindex{55.5}{11.5}{1}
\zindex{44.5}{4.5}{l}
\zindex{55.5}{4.5}{1}
\zindex{38}{11.5}{b_2}
\zindex{38.5}{4.5}{b_1}
\zindex{62.5}{11.5}{a_1}
\zindex{62.5}{4.5}{a_2}

\zindex{90.2}{7.5}{l}
\zindex{93.8}{7.5}{k}
\zindex{91.5}{13.5}{k}
\zindex{102.5}{13.5}{1}
\zindex{84}{11.5}{b_1}
\zindex{86}{14.8}{b_2}
\zindex{109.5}{13.5}{a_1}
\zindex{109.5}{6.5}{a_2}

\zindex{8}{9.8}{l}
\zindex{8}{6.2}{k}
\zindex{14}{14.8}{la_1}
\zindex{17}{12}{la_2}
\zindex{13.2}{1.5}{b_2}
\zindex{16.2}{4}{b_1}

\end{pict}%
The same $G$-tree results if we first perform slides and then expand,
including $\partial_0(e_1)$ in the same side as $\partial_1(e_1)$, and
then exchange the names of $e$ and $e_1$. 
\end{case1}
\begin{case2} The edge $e_1$ has distinct endpoints and is incident to
only one endpoint of $e$. 
Let $\{f_i\}$ be the edges with $\partial_0(f_i) = \partial_0(e)$ just
before the slide move (not including $e_1$). We replace the 
expansion and slide by slides and an expansion as follows: first slide
each $f_i$ over $e_1$, then expand at $\partial_1(e_1)$ so that the
new expansion edge $e$ separates $\{f_i\}$ from the rest of the edges at
$\partial_1(e_1)$. For example: 
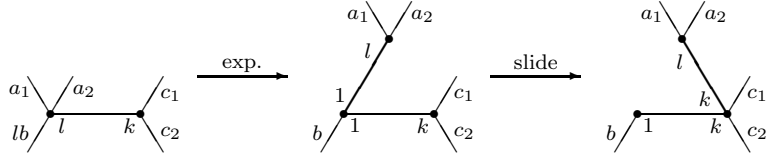
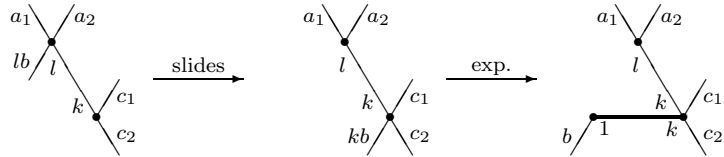
\begin{pict}{100}{20}
\thicklines
\put(5,5){\circle*{1}}
\put(17,5){\circle*{1}}

\put(44,5){\circle*{1}}
\put(50,15){\circle*{1}}
\put(56,5){\circle*{1}}
\put(44,5){\line(3,5){6}}

\put(83,5){\circle*{1}}
\put(89,15){\circle*{1}}
\put(95,5){\circle*{1}}
\put(95,5){\line(-3,5){6}}

\thinlines
\put(5,5){\line(1,0){12}}
\put(5,5){\line(-3,-5){3}}
\put(5,5){\line(-3,5){3}}
\put(5,5){\line(3,5){3}}
\put(17,5){\line(3,-5){3}}
\put(17,5){\line(3,5){3}}

\zindex{1}{2.5}{lb}
\zindex{1}{8}{a_1}
\zindex{9.5}{8}{a_2}
\zindex{6.5}{3.5}{l}
\zindex{15.5}{3.5}{k}
\zindex{21}{7.5}{c_1}
\zindex{21}{2.5}{c_2}

\put(24.5,10){\vector(1,0){12}}
\zindex{30.5}{11.5}{\mbox{exp.}}

\put(44,5){\line(1,0){12}}
\put(44,5){\line(-3,-5){3}}
\put(50,15){\line(-3,5){3}}
\put(50,15){\line(3,5){3}}
\put(56,5){\line(3,-5){3}}
\put(56,5){\line(3,5){3}}

\zindex{40.5}{2.5}{b}
\zindex{46}{18}{a_1}
\zindex{54.5}{18}{a_2}
\zindex{47.2}{13.5}{l}
\zindex{43.5}{7.5}{1}
\zindex{45.5}{3.5}{1}
\zindex{54.5}{3.5}{k}
\zindex{60}{7.5}{c_1}
\zindex{60}{2.5}{c_2}

\put(63.5,10){\vector(1,0){12}}
\zindex{69.5}{12}{\mbox{slide}}

\put(83,5){\line(1,0){12}}
\put(83,5){\line(-3,-5){3}}
\put(89,15){\line(-3,5){3}}
\put(89,15){\line(3,5){3}}
\put(95,5){\line(3,-5){3}}
\put(95,5){\line(3,5){3}}

\zindex{79.5}{2.5}{b}
\zindex{85}{18}{a_1}
\zindex{93.5}{18}{a_2}
\zindex{88.7}{12}{l}
\zindex{92}{7}{k}
\zindex{84.5}{3.5}{1}
\zindex{93.5}{3.5}{k}
\zindex{99}{7.5}{c_1}
\zindex{99}{2.5}{c_2}

\end{pict}%
becomes
\begin{pict}{100}{20}
\thicklines
\put(11,15){\circle*{1}}
\put(17,5){\circle*{1}}

\put(50,15){\circle*{1}}
\put(56,5){\circle*{1}}

\put(83,5){\circle*{1}}
\put(89,15){\circle*{1}}
\put(95,5){\circle*{1}}
\put(83,5){\line(1,0){12}}

\thinlines
\put(17,5){\line(-3,5){6}}
\put(11,15){\line(-3,-5){3}}
\put(11,15){\line(-3,5){3}}
\put(11,15){\line(3,5){3}}
\put(17,5){\line(3,-5){3}}
\put(17,5){\line(3,5){3}}

\zindex{7}{12.5}{lb}
\zindex{7}{18}{a_1}
\zindex{15.5}{18}{a_2}
\zindex{11.2}{11.5}{l}
\zindex{14.5}{6.5}{k}
\zindex{21}{7.5}{c_1}
\zindex{21}{2.5}{c_2}

\put(24.5,10){\vector(1,0){12}}
\zindex{30.5}{12}{\mbox{slides}}

\put(56,5){\line(-3,5){6}}
\put(56,5){\line(-3,-5){3}}
\put(50,15){\line(-3,5){3}}
\put(50,15){\line(3,5){3}}
\put(56,5){\line(3,-5){3}}
\put(56,5){\line(3,5){3}}

\zindex{51.8}{2.5}{kb}
\zindex{46}{18}{a_1}
\zindex{54.5}{18}{a_2}
\zindex{49.7}{12}{l}
\zindex{53}{7}{k}
\zindex{60}{7.5}{c_1}
\zindex{60}{2.5}{c_2}

\put(63.5,10){\vector(1,0){12}}
\zindex{69.5}{11.5}{\mbox{exp.}}

\put(95,5){\line(-3,5){6}}
\put(83,5){\line(-3,-5){3}}
\put(89,15){\line(-3,5){3}}
\put(89,15){\line(3,5){3}}
\put(95,5){\line(3,-5){3}}
\put(95,5){\line(3,5){3}}

\zindex{79.5}{2.5}{b}
\zindex{85}{18}{a_1}
\zindex{93.5}{18}{a_2}
\zindex{88.7}{12}{l}
\zindex{92}{7}{k}
\zindex{84.5}{3.5}{1}
\zindex{93.5}{3.5}{k}
\zindex{99}{7.5}{c_1}
\zindex{99}{2.5}{c_2}

\end{pict}%
As in Case 1, the names of $e$ and $e_1$ must be exchanged after the new
moves. This procedure works whether $e_1$ is incident to $\partial_0(e)$
or to $\partial_1(e)$. 
\end{case2}
\begin{case3} The edge $e_1$ is a loop incident to $\partial_1(e)$. Then
the procedure from Case 2 works. The two moves are replaced by a
sequence of slides (around the loop $e_1$) followed by an expansion. 
\end{case3}
\begin{case4} The edge $e_1$ is a loop incident to $\partial_0(e)$. Let
$l = i(\overline{e})$. Since $\partial_0(e)$ is the end of $e$ that
slides over $e_1$ and $i(e)=1$, the loop $e_1$ must be an ascending
loop (before and after the slide). Note that before the expansion of
$e$, the indices of $e_1$ were $l$ times their current values; hence
$e_1$ was originally a virtually ascending loop. Let $k$ be the modulus
of $e_1$. If $k \not= \pm 1$ then alternative (i) holds. If $k=1$ then
the slide move may simply be omitted. If $k=-1$ then first perform the
slide moves as described in Case 2, and then expand the edge $e$ as
before, but with $i(\overline{e}) = l$ and $i(e) = -1$.
\qed \end{case4}
\end{proof2}

\begin{proposition}\label{expcollapse}
Suppose an expansion creating the edge $e$ is followed by the collapse of
an edge $e'$. Then either
\begin{enum}
\item[\textup{(i)}] $e'$ is a strict ascending loop before the expansion
  and $e$ is a strict ascending loop after the collapse, 
\item[\textup{(ii)}] both moves may be deleted, 
\item[\textup{(iii)}] both moves may be replaced by a sequence of slides,
or 
\item[\textup{(iv)}] the collapse may be performed before the expansion
move. 
\end{enum}
\end{proposition}

\begin{proof2}
If $e = e'$ or $e = \overline{e}'$ then clearly (ii) holds. Otherwise $e$
and $e'$ are distinct, proper edges just after the expansion and before
the collapse. If they do not meet then conclusion (iv) holds. 

Now assume that $e$ and $e'$ meet in one or two vertices. Orient both
edges so that $i(e) = i(e') = 1$. 
\begin{case1} The edges $e$ and $e'$ have two vertices in common. If
$\partial_0(e) = \partial_1(e')$ and $\partial_1(e) = \partial_0(e')$
then alternative (i) holds. Otherwise, if $\partial_0(e) =
\partial_0(e')$ and $\partial_1(e) = \partial_1(e')$ then set $k =
i(\overline{e})$ and $l = (\overline{e}')$. The moves have the form: 
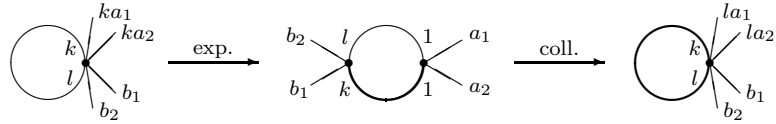
\begin{pict}{102}{16}
\scriptsize
\thicklines
\put(10,8){\circle*{1}}
\put(45,8){\circle*{1}}
\put(55,8){\circle*{1}}
\put(50,8){\oval(10,10)[b]}
\put(88,8){\circle{10}}

\put(93,8){\circle*{1}}

\thinlines
\put(5,8){\circle{10}}
\put(50,8){\circle{10}}

\put(21,8){\vector(1,0){12}}
\put(67,8){\vector(1,0){12}}
\zindex{27}{9.5}{\mbox{exp.}}
\zindex{73}{10}{\mbox{coll.}}
\put(10,8){\line(1,6){1}}
\put(10,8){\line(1,1){4}}
\put(10,8){\line(1,-1){4}}
\put(10,8){\line(1,-6){1}}

\put(45,8){\line(-5,3){5}}
\put(45,8){\line(-5,-3){5}}
\put(55,8){\line(5,3){5}}
\put(55,8){\line(5,-3){5}}

\put(93,8){\line(1,6){1}}
\put(93,8){\line(1,1){4}}
\put(93,8){\line(1,-1){4}}
\put(93,8){\line(1,-6){1}}

\zindex{44.5}{11.5}{l}
\zindex{55.5}{11.5}{1}
\zindex{44.5}{4.5}{k}
\zindex{55.5}{4.5}{1}
\zindex{38}{11.5}{b_2}
\zindex{38.5}{4.5}{b_1}
\zindex{62.5}{11.5}{a_1}
\zindex{62.5}{4.5}{a_2}

\zindex{8}{10}{k}
\zindex{8}{6}{l}
\zindex{14}{14.8}{ka_1}
\zindex{17}{12}{ka_2}
\zindex{13.2}{1.5}{b_2}
\zindex{16.2}{4}{b_1}

%% \zindex{91}{10}{k}
%% \zindex{91}{6}{l}
%% \zindex{96}{14.8}{b_1}
%% \zindex{99}{12}{b_2}
%% \zindex{96.7}{1.5}{la_2}
%% \zindex{100}{4}{la_1}
\zindex{91}{10}{k}
\zindex{91}{6}{l}
\zindex{96.7}{14.8}{la_1}
\zindex{99.8}{12}{la_2}
\zindex{96.2}{1.5}{b_2}
\zindex{99.2}{4}{b_1}

\end{pict}%
Evidently the moves may be replaced by slides around the loop $e'$. 
\end{case1}
\begin{case2} The edges $e$ and $e'$ meet in one vertex. Again let $k =
i(\overline{e})$ and $l = (\overline{e}')$. There are four
configurations. If $\partial_0(e) = \partial_0(e')$ then we see:
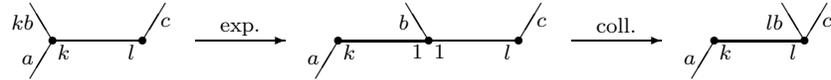
\nopagebreak 
\begin{pict}{110}{10}
\scriptsize
\thicklines
\put(6,5){\circle*{1}}
\put(18,5){\circle*{1}}

\put(44,5){\circle*{1}}
\put(56,5){\circle*{1}}
\put(68,5){\circle*{1}}
\put(44,5){\line(1,0){12}}

\put(94,5){\circle*{1}}
\put(106,5){\circle*{1}}
\put(94,5){\line(1,0){12}}

\thinlines
\put(6,5){\line(1,0){12}}
\put(56,5){\line(1,0){12}}

\put(25,5){\vector(1,0){12}}
\put(75,5){\vector(1,0){12}}
\zindex{31}{6.5}{\mbox{exp.}}
\zindex{81}{7}{\mbox{coll.}}

\put(6,5){\line(-3,5){3}}
\put(6,5){\line(-3,-5){3}}
\put(18,5){\line(3,5){3}}
\zindex{7.5}{3.5}{k}
\zindex{16.5}{3.5}{l}
\zindex{2.7}{2.5}{a}
\zindex{2}{7.5}{kb}
\zindex{21}{7.5}{c}

\put(56,5){\line(-3,5){3}}
\put(44,5){\line(-3,-5){3}}
\put(68,5){\line(3,5){3}}
\zindex{45.5}{3.5}{k}
\zindex{54.5}{3.5}{1}
\zindex{57.5}{3.5}{1}
\zindex{66.5}{3.5}{l}
\zindex{40.7}{2.5}{a}
\zindex{52.7}{7.5}{b}
\zindex{71}{7.5}{c}

\put(106,5){\line(-3,5){3}}
\put(94,5){\line(-3,-5){3}}
\put(106,5){\line(3,5){3}}
\zindex{95.5}{3.5}{k}
\zindex{104.5}{3.5}{l}
\zindex{90.7}{2.5}{a}
\zindex{102}{7.5}{lb}
\zindex{109}{7.5}{c}

\end{pict}%
We may replace the two moves by slides over the edge $e'$. 

In the other three configurations the collapse may be performed before
the expansion. To illustrate, the case $\partial_1(e) = \partial_1(e')$
has the following configuration: \nopagebreak
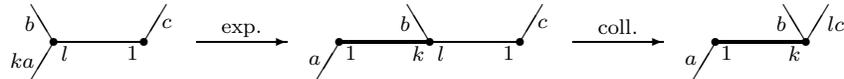
\begin{pict}{111}{10}
\scriptsize
\thicklines
\put(6,5){\circle*{1}}
\put(18,5){\circle*{1}}

\put(44,5){\circle*{1}}
\put(56,5){\circle*{1}}
\put(68,5){\circle*{1}}
\put(44,5){\line(1,0){12}}

\put(94,5){\circle*{1}}
\put(106,5){\circle*{1}}
\put(94,5){\line(1,0){12}}

\thinlines
\put(6,5){\line(1,0){12}}
\put(56,5){\line(1,0){12}}

\put(25,5){\vector(1,0){12}}
\put(75,5){\vector(1,0){12}}
\zindex{31}{6.5}{\mbox{exp.}}
\zindex{81}{7}{\mbox{coll.}}

\put(6,5){\line(-3,5){3}}
\put(6,5){\line(-3,-5){3}}
\put(18,5){\line(3,5){3}}
\zindex{7.5}{3.5}{l}
\zindex{16.5}{3.5}{1}
\zindex{1.8}{2.5}{ka}
\zindex{2.8}{7.5}{b}
\zindex{21}{7.5}{c}

\put(56,5){\line(-3,5){3}}
\put(44,5){\line(-3,-5){3}}
\put(68,5){\line(3,5){3}}
\zindex{45.5}{3.5}{1}
\zindex{54.5}{3.5}{k}
\zindex{57.5}{3.5}{l}
\zindex{66.5}{3.5}{1}
\zindex{40.7}{2.5}{a}
\zindex{52.7}{7.5}{b}
\zindex{71}{7.5}{c}

\put(106,5){\line(-3,5){3}}
\put(94,5){\line(-3,-5){3}}
\put(106,5){\line(3,5){3}}
\zindex{95.5}{3.5}{1}
\zindex{104.5}{3.5}{k}
\zindex{90.7}{2.5}{a}
\zindex{102.7}{7.5}{b}
\zindex{110}{7.5}{lc}

\end{pict}%
and it is easy to see that the collapse may be performed first. The
remaining two cases are entirely similar. 
\qed \end{case2}
\end{proof2}

\begin{theorem}\label{slidethm} 
Let $X$ and $Y$ be reduced non-elementary cocompact generalized
Baumslag--Solitar trees with group $G$, and suppose that $q({G}) \cap \Z
= 1$. Then $X$ and $Y$ are related by slide moves. 
\end{theorem}

\begin{proof}
The property $q({G}) \cap \Z = 1$ guarantees that no
generalized Baumslag--Solitar decomposition of ${G}$ contains strict
virtually ascending loops. Starting with a sequence of moves from $X$ to
$Y$ (given by Theorem \ref{defthm}) we claim that Propositions
\ref{expslide} and \ref{expcollapse} can be applied to obtain a new
sequence consisting of collapses, followed by slides, followed by
expansions. To see this, note that case (i) of either proposition cannot
occur. Therefore expansions can be pushed forward past slides (by
\ref{expslide}) and past collapses (by \ref{expcollapse}), and collapses
can be pulled back before slides (by \ref{expslide} applied to the
reverse of the sequence of moves). That is, we have the replacement rules
$ES \to S^* E$, $EC \to (S^* \mbox{ or } CE)$, and $SC \to CS^*$, where
$E$ and $C$ denote expansion and collapse moves respectively and $S^*$
denotes a (possibly empty) sequence of slide moves. 

The algorithm for simplifying a sequence of moves is to repeatedly
perform either of the following two steps, until neither applies. The
first step is to find the first collapse move that is preceded by an
expansion or slide, and apply the replacement $EC \to (S^* \mbox{ or }
CE)$ or $SC \to CS^*$ accordingly. The second step is to find the last
expansion move that is followed by a collapse or slide and apply the
replacement $EC \to (S^* \mbox{ or } CE)$ or $ES \to S^* E$. This
procedure terminates, in a sequence of the form $C^* S^* E^*$. Then since
$X$ and $Y$ are reduced, the new sequence of moves has no collapses or
expansions. 
\end{proof}

\section{The isomorphism problem}\label{sec8}

Next we approach the problem of classifying generalized Baumslag--Solitar
groups. At the minimum, a classification should include an algorithm for
determining when two indexed graphs define the same group. This is the
problem considered here. 

For certain generalized Baumslag--Solitar groups the isomorphism problem
is trivial. This occurs when the deformation space contains only one
reduced tree (such a tree is called \emph{rigid}). The basic rigidity
theorem for generalized Baumslag--Solitar trees was proved independently
in \cite{pettet,gilbertetal,forester:trees} and it states that
$(A,i)_{\Z}$ is rigid if there are no divisibility relations at any
vertex. Levitt \cite{levitt:char} has extended this result by giving a
complete characterization of trees that are rigid. Then to solve the
isomorphism problem for such groups one simply makes the trees reduced
and compares them directly (cf. Remark \ref{signs}). 

In this section we solve the isomorphism problem for the case of
generalized Baumslag--Solitar groups having no non-trivial integral
moduli. The general case is still open. 

\begin{lemma}\label{indexbound} 
Let $Q \subset \Q^{\times}_{>0}$ be a finitely generated subgroup such
that $Q \cap \Z = 1$. Then for any $r \in \Q$ the set $rQ \cap \Z$
is finite. 
\end{lemma}

\begin{proof} 
We consider $\Q^{\times}_{>0}$ as a free $\Z$-module with basis the
prime numbers, via prime decompositions. Note that a positive rational
number is an integer if and only if it has nonnegative coordinates in
$\Q^{\times}_{>0} = \Z \oplus \Z \oplus \cdots$, and so the positive
integers comprise the first ``orthant'' of $\Q^{\times}_{>0}$. 

We are given that $Q$ meets the first orthant only at the origin. By 
taking tensor products with $\R$ we may think of $\Q^{\times}_{>0}$ as 
a vector space and $Q$ a finite dimensional subspace. Since
multiplication by $r$ is a translation in $\Z \oplus \Z \oplus \cdots$,
we have that $rQ$ is an affine subspace parallel to $Q$. It suffices to
show that this affine subspace meets the first orthant in a compact set. 

This is clear if $Q$ is a codimension $1$ subspace of a coordinate
subspace $\R \oplus \cdots \oplus \R$, because the subspace would have a
strictly positive normal vector, and then $rQ$ would meet the first
orthant of $\R \oplus \cdots \oplus \R$ in a simplex or a point (or not at
all). Otherwise we can choose a coordinate subspace $\R \oplus \cdots
\oplus \R$ containing $Q$ and then enlarge $Q$ to make it codimension
$1$, preserving the property that it meets the first orthant only at the
origin. The result then follows easily. 
\end{proof}

\begin{theorem}\label{finiteness}
Let $G$ be a finitely generated generalized Baumslag--Solitar group. If
$q(G) \cap \Z = 1$ then there are only finitely many reduced graphs of
groups $(A,i)_{\Z}$ with fundamental group $G$. 
\end{theorem}

\begin{proof}
If $G$ is elementary then there are only four reduced graphs of groups 
whose universal covering trees have at most two ends and the result is 
clear. These are: a single vertex, a loop with indices $\pm 1$ (two
cases: equal signs or opposite signs), and an interval with indices $\pm
2$. 

If $G$ is non-elementary then any two reduced trees are
related by slide moves, by Theorem \ref{slidethm}. In particular there
are only finitely many possible quotient graphs. Thus we may consider
sequences of slide moves in which every edge returns to its original
position in the quotient graph. To prove the theorem it then suffices to
show that after such a sequence, there are only finitely many possible
values for each edge index $i(e)$. 

Suppose $X$ is a generalized Baumslag--Solitar tree with group $G$ and
$e$ is an edge of $X$ with initial vertex $v$. Consider a sequence of
slide moves after which $e$ has initial vertex $gv$ for some $g \in G$.
Then we have $G_e \subset (G_v \cap G_{gv})$. This implies that 
\begin{equation*}
\begin{split}
[G_v:G_e] \ &= \ [G_v:(G_v \cap
G_{gv})] \, [(G_v \cap G_{gv}):G_e] \\
&= \ q(g)\, [G_{gv}:G_e], 
\end{split}
\end{equation*}
or equivalently $[G_{gv}:G_e] = [G_v:G_e] \, q(g^{-1})$.
Therefore, since $i(e) = \pm [G_v:G_e]$ before the slide moves, the new
index after the moves is an element of the set $\pm i(e) q(G) \cap
\Z$. This set is finite by Lemma \ref{indexbound}. 
\end{proof}

\begin{corollary}\label{isomproblem}
There is an algorithm which, given finite graphs of groups $(A,i)_{\Z}$ and
$(B,j)_{\Z}$ such that $(A,i)_{\Z}$ has no non-trivial integral moduli,
determines whether the associated generalized Baumslag--Solitar groups
are isomorphic. 
\end{corollary}

\begin{proof}
First make $(A,i)_{\Z}$ and $(B,j)_{\Z}$ reduced. If one or both is
elementary then it is a simple matter to check for isomorphism. Among the
reduced elementary graphs of groups, a single vertex has group $\Z$, a
loop with both indices equal to $1$ has group $\Z \times \Z$, and the
remaining two cases yield the Klein bottle group. 

Otherwise, by Theorem \ref{slidethm}, the groups are isomorphic if and
only if there is a sequence of slide moves taking $(A,i)_{\Z}$ to
$(B,j)_{\Z}$. Now consider the set of graphs of groups related to
$(A,i)_{\Z}$ by slide moves. This is the vertex set of a connected graph
$\Gs$ whose edges correspond to slide moves. We claim that every vertex
of $\Gs$ is a reduced graph of groups. This then implies that $\Gs$ is
finite by Theorem \ref{finiteness}. To prove the claim one observes,
using Proposition \ref{gbsmoves}, that an edge in a reduced graph of
groups cannot be made collapsible during a slide move unless it slides
over a strict ascending loop, but there are no such loops because the
group has no non-trivial integral moduli. 

Now search $\Gs$ by performing all possible sequences of slide moves of
length $n$, for increasing $n$ until no new graphs of groups are
obtained. Then the two generalized Baumslag--Solitar groups are
isomorphic if and only if $(B,j)_{\Z}$ has been found by this point. 
\end{proof}

\end{document}